# Dynamic Workforce Scheduling and Relocation in Hyperconnected Parcel Logistic Hubs

**Yujia Xu, Yiguo Liu, and Benoit Montreuil**
*Physical Internet Center, ISyE, Georgia Institute of Technology*
*Atlanta, Georgia, 30332, USA*

## Abstract

With the development of e-commerce during the Covid-19 pandemic, one of the major challenges for many parcel logistics companies is to design reliable and flexible scheduling algorithms to meet uncertainties of parcel arrivals as well as manpower supplies in logistic hubs, especially for those depending on workforce greatly. Currently, most labor scheduling is periodic and limited to single facility, thus the number of required workers in each hub is constrained to meet the peak demand with high variance. We approach this challenge, recognizing that not only workforce schedules but also working locations could be dynamically optimized by developing a dynamic workforce scheduling and relocation system, fed from updated data with sensors and dynamically updated hub arrival demand predictions. In this paper, we propose novel reactive scheduling heuristics to dynamically match predicted arrivals with shifts at hyperconnected parcel logistics hubs. Dynamic scheduling and allocation mechanisms are carried out dynamically during delivery periods to spatiotemporally adjust the available workforce. We also include penalty costs to keep parcels sorted in time and scheduling adjustments are made in advance to allow sufficient time for crew planning. To assess the proposed methods, we conduct comprehensive case studies based on real-world parcel logistic networks of a logistic company in China. The results show that our proposed approach can significantly outperform traditional workforce scheduling strategies in hubs with limited computation time.

## Keywords

## 1. Introduction

Shift scheduling is one of the oldest operation research problems and plays an important role in manufacturing and service industries. The goal of the problem is to ensure efficient use of resources, to achieve balanced workload distribution and to meet individual needs as much as possible [1]. The shift scheduling problem can be divided into two broad categories based on the type of workload they consider: task-coverage problems and workload-coverage problems. For task-coverage problems, tasks to be executed during the planning horizon are known and must be assigned to shifts or workers. In contrast, in workload-coverage problems, tasks to be executed are unknown and the demand for employees is forecasted based on expected workloads and workers are scheduled to cover these predicted personnel demands [2].

As traditional supply chains are increasingly becoming intelligent with more objects embedded with sensors and better communication, intelligent decision making and automation capabilities, the new smart supply chains present unprecedented opportunities for achieving cost reduction and enhancing efficiency improvement [3]. In the past few decades, dynamic scheduling problems have attracted widespread attention, focused on making timely decisions considering real-time system status with uncertainties. Dynamic scheduling problems can be categorized into three types for dealing with uncertainty in schedule disruptions. Completely reactive scheduling makes local decisions in real-time using heuristic rules to respond quickly to changes of workshop state and provide guidance for the next operation of actual production system [4]. Predictive-reactive scheduling generates a preliminary predictive schedule prior to execution and may require modification following a schedule disruption [5]. Robust scheduling generates a prescheduling scheme that can accommodate various possible disruptions based on forecasting on future dynamic events, to make sure the performance of realized scheduling system do not deteriorate with dynamic events [6].

The Physical Internet (PI) was introduced by Montreuil [7] as 'an open global logistics system founded on physical, digital and operational interconnectivity through encapsulation, interfaces and protocols', and thus defines a new opportunity for supply chain design and operations, enabling seamless open asset sharing and flow consolidation [8]. In the context of urban logistics, it is materialized through a multi-tier urban logistics web [9] whose nodes are hyperconnected logistic hubs enabling the sorting, consolidation, transshipment, and crossdocking of goods [10]. The multi-tier structure has access hubs interconnect unit zones, local hubs interconnect local cells, and gateway hubs interconnect urban areas [10], which provides an opportunity for moving workers across these hubs. Some of the hubs are close to each other but on different planes in the meshed networks, which means they have different arrival or departure patterns and workforce demand peaks during a day. For instance, arrival time of parcels at access hubs depends on when customers require them to be picked up and the arrival time of parcels at gateway hubs depends on the schedules of transport trucks, subways, or trains.

In this paper, we propose a dynamic workforce scheduling and relocation system to match predicted demand with shifts in real-time using heuristic rules to respond quickly to changes of predicted arrivals at hubs and provide guidance for centralized workforce assignment. Shifts are scheduled to cover the predicted workloads at hubs and online scheduling problems are solved dynamically based on real-time status including the assigned shifts as well as future tasks. We also rely on a rolling horizon approach to address the presence of uncertainty, which is to implement a reactive scheduling method that iteratively solves the deterministic problem by moving forward the optimization horizon in every iteration; assuming that the status of the system is updated as soon as the different uncertain or not accurate enough parameters become known, the schedule can be optimized for the new resulting scenario [11]. Although some related works exist, for example, [12] proposed demand prediction and workforce allocation models to improve airport screening operations and [13] developed workforce scheduling and routing models for service providers, we believe that this is the first work that shows the feasibility, efficiency and reliability of the proposed dynamic workforce scheduling and relocation system for hyperconnected logistics hubs considering parcels' dynamically predicted arrival time and maximum dwell time in real-world logistics networks.

## 2. Problem Description and Methodology

We assume that workload prediction for hubs in the hyperconnected logistic networks are provided, updated every $\epsilon$ minutes. Since the exact capabilities for solving the scheduling problem for large logistics networks which contains hundreds of hubs are limited, we present a reactive scheduling model in a rolling horizon to iteratively generate continuous shifts and shift combinations. In addition, unlike the traditional scheduling problem for a single logistic hub, we incorporate the possibility of labor sequentially working at multiple nearby hubs within daily shifts. The flowchart of the developed methodology can be found in **Figure 1**. The overall objective is to cover workload at hubs with minimized cost and penalties including parcels' lateness and emergency hiring.

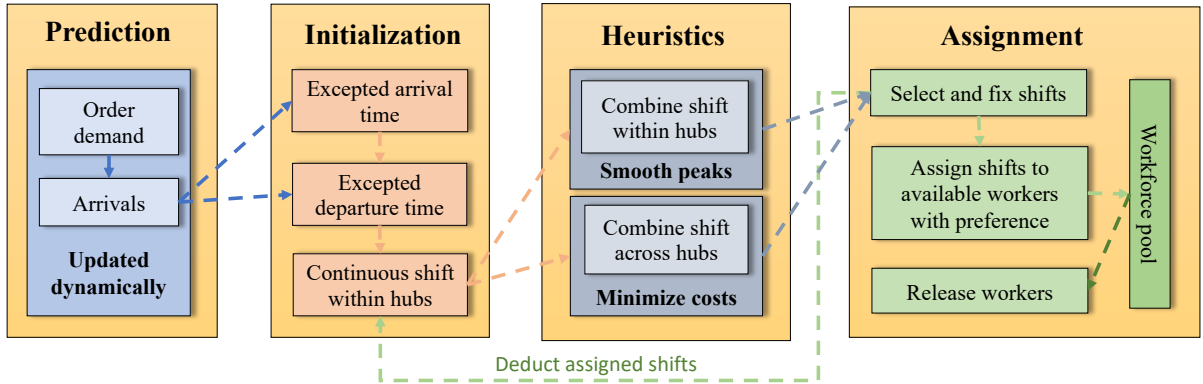

**Figure 1:** Flowchart of the developed reactive scheduling methodology

### 2.1 Heuristic Approach

Given the dynamically updated workload prediction during the planning horizon, shifts with the maximum possible length for each hub are initialized accordingly. We assume that after a parcel arrived at a hub, there is a specified maximum dwell time $\zeta$, during which it should be processed and then leave the hub. Thus, when we combine the short shifts into a long shift, workforce demand peaks can be smoothed utilizing this dwell time. In other words, parcels' arrival peaks do not always result in workforce peaks since parcels are not necessary to be sorted and processed right after they arrive, and the delayed work may lead to smoothed workforce demand. After combing possible shifts within hubs, we use greedy algorithms to find and combine short shifts across nearby hubs, during which we also take the transportation modes into consideration. Shifts are combined between hubs greedily until the distance-related travelling cost is larger than hiring a new worker. Once a continuous shift or shift combination is chosen, the centralized workforce agent assigns an available worker to it. The procedure of this reactive scheduling model is detailed below in pseudocode format.

```
PROCEDURE
/* Part I: Initialize the labor demand in each time unit by dividing the
number of arrived parcels by labor efficiency, x = [x_0, x_1 ... x_n]. */
shifts = [], start = 0
WHILE start < n THEN
    IF x[start] = 0 THEN
        start = start + 1
    x[start] = x[start] – 1
    end = start + 1
    /* find the shift with the maximum possible length*/
    WHILE (end < n) and (x[end] > 0) and (end – start <=
    maximum working hour ρ) THEN
        x[end] = x[end] – 1
        end = end + 1
    shifts = shifts ∪ {[start, end]}
RETURN shifts

/* Part II: combine the workload with minimized number of shifts
and maximized shift total length utilizing the parcels' dwell time at
hubs to smooth the labor demand peaks */
Combined = [], unit labor demand x = [x_0, x_1 ... x_n].
FOR shift in shifts
    IF shift[end] – shift[start] = maximum working hour ρ
THEN
        Combined. append(shift)
        x[shift[start] : shift[end]] –=1
start = 0
WHILE start < n THEN
    IF x[start] = 0 THEN
        start += 1
        Continue
    end = start + 1

    WHILE (end < n) AND (end – start < maximum working
    hour ρ) THEN
        Find next non-zero hour and add to shift
        IF x[end] = 0 AND x[end-1] != 0 THEN
            /*smooth labor demand by moving
count to next time unit */
            Move count from end-1 to end
    Combined = Combined ∪ {shifts}
    x[shift[start] : shift[end]] –= 1
RETURN Combined

/* Part III: greedily merge the shifts across hubs */
Combined = []
H = {[hub pairs]}/* all possible moving pairs with ascending
distance but smaller than largest moving distance δ. */
FOR hub pair in H
    FOR shifts in the pair
        IF (shift hours do not overlap) AND (moving
time <= shifts gap <= maximum gap) AND (total working time <=
maximum working hour ρ) THEN
            Combined = Combined ∪ {shifts}
            Remove shifts from hubs
RETURN Combined

/* Part IV: select shifts and assign them by to available candidates
from the pool; release workers if finished*/
For shifts in Combined:
    IF Value(shift) >= threshold γ THEN /* Value(shifts) is
calculated by value function defined */
        assign an available worker in the pool to it according to
preference if possible
```

## 2.2 Rolling Horizon Approach

A rolling horizon approach to address the presence of uncertainty then follows, seeking to assign cheap and efficient shifts to workers iteratively by moving forward the planning horizon with updated demand prediction at every iteration. Then the selection of shifts in each planning horizon becomes important. On one hand, fixing shifts too early may lead to overstaffing in presence of significant demand prediction intervals. On the other hand, selecting shifts too late may result in high emergency penalty and leave no time for crew planning. Thus, we evaluate the shifts using the following value function:

$$\begin{cases} Value = \alpha \times \frac{\text{Duration threshold } \tau \text{ to fix a shift}}{\text{Shift start time } - \text{ Current time } t} + \beta \times \frac{\text{Working time}}{\text{Maximum working time } \rho} + \gamma \times \frac{\text{Working time}}{\text{Resting time}} \\ \alpha + \beta + \gamma = 1 \end{cases}$$

If the value of a continuous shift or a shift combination is more than a threshold $\delta$ ($0 \leq \delta \leq 1$), it would be fixed and assigned to an available worker from the workforce pool. The parameters $\alpha$, $\beta$, $\gamma$, $\delta$ and $\tau$ can be tuned to minimize the total costs according to different prediction scenarios. In a rolling horizon, every time workload predictions are updated and some shifts are selected and assigned to a worker, the corresponding workload is deducted from new prediction results, which become input for the next planning horizon. In addition, once a worker has finished her/his assigned shifts and is available for her/his next day shifts, he/she would be released to the workforce pool.

## 3. Experimental Results

To assess the benefits of our proposed model, we conduct experiments using logistics networks from a logistics company in China, leveraging an urban logistics simulator to dynamically collect parcels' arrivals at hyperconnected

*Xu, Liu, and Montreuil*

hubs [10]. In our experiments, 52 local hubs and gateway hubs from the company's logistics networks in Shenzhen, China are included. We set the planning horizon to be 24 hours and there exist 1,173,253 arrivals at hubs during the first day. We assume the maximum dwell time $\zeta = 1$ hour, the maximum working time $\rho = 8$ hours, the working efficiency $\mu = 150$ parcels/hour, and our models to be run every $\epsilon = 60$ minutes utilizing updated predictive results. It takes about a minute to run the methodology for a whole day with the rolling horizons.

Since the prediction of parcels' arrival time at hubs can be dynamically improved as they approach to the hubs, in the experiment we assume that the prediction of number of arrivals $p$ of time $t_1$ conducted at time $t_0$ can be written as

$$p = max\left(0, \frac{C}{100} * (u * (t_1 - t_0) + 100)\right)$$

where $u \sim Uniform[-1,1]$ and $C$ is the actual number of arrivals at time $t_1$. The parameters to fix combined shifts during the rolling horizon are tuned like α = 0.4, β = 0.3, γ = 0.3, δ = 0.9 and τ = 4. As these influence performance, different settings may be further tested.

### 3.1 Complexity of Assigned Shifts

Results show that if we restrict shifts to be within a single hub, 24.2% of shifts assigned have multiple durations, which means workers need to take a rest within a hub for some time and then continue to work, while 75.8% of the shifts are continued without resting. By allowing movements across hubs, we assign 9.1% of total shifts for workers moving across nearby hubs, as shown in **Figure 2**. In addition, even though most of the shifts assigned within hubs have maximum allowed duration of 8 hours, allowing movements across hubs improve the average working duration by combining short shifts into a long shift. The comparison of the number of hours in daily of shifts between with and without movements across hubs is also shown in **Figure 2**.

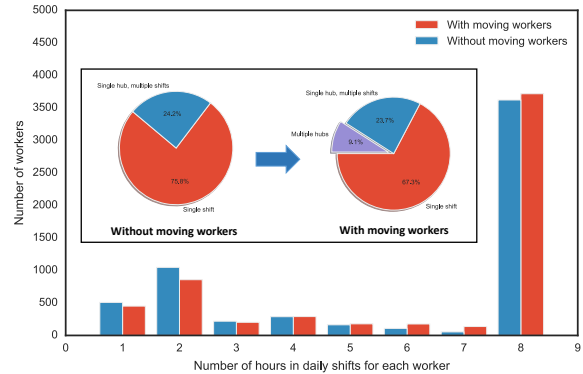

**Figure 2:** Complexity of assigned shifts

### 3.2 Arrivals and Workers Assigned at Hubs

Shifts are dynamically assigned according to workload prediction at hubs. To show the actual number of parcels arrived as well as working and resting workers at hubs during a day, we present their dynamic evolution in two gateway hubs and two local hubs in a simulation instance. As shown in **Figure 3**, the black lines indicate the actual number of parcels arrived at hubs while the purple lines show the total number of workers at hubs. Among them, some are working while some are taking rests, shown as red and blue dashed lines accordingly. It is shown that the number of workers at hubs shares the same pattern with arrivals but the arrival peaks are smoothed. Also, workers tend to take rest around arrival valleys.

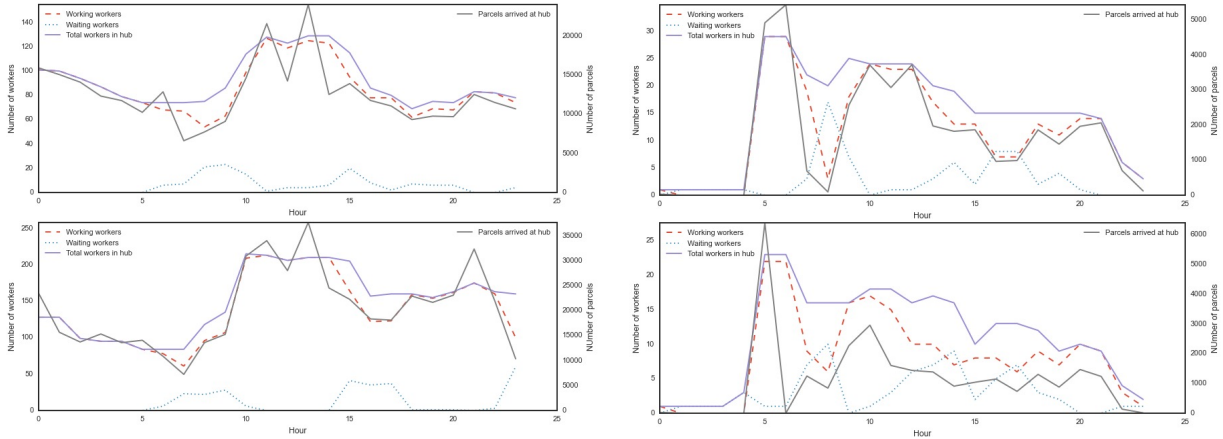

**Figure 3:** Parcels arrived and workers assigned at two larger gateway hubs (left) and two smaller local hubs (right)

## 3.2 Workers in Hubs and Flow across Hubs

We also plot the locations of hubs as well as how workers appear and move across hubs during a day. As shown in **Figure 4**, the size of the circles indicates the number of shifts assigned in each hub while the width of lines signifies the number of workers transported from one hub to another hub every 6 hours. The three large circles in each picture mean that many shifts are assigned to the three gateway hubs to cover comparatively large quantities of arrivals. Also, many workers move between gateway hubs and local hubs from 12 AM to 12 PM because of their different parcel arrival patterns during this time span.

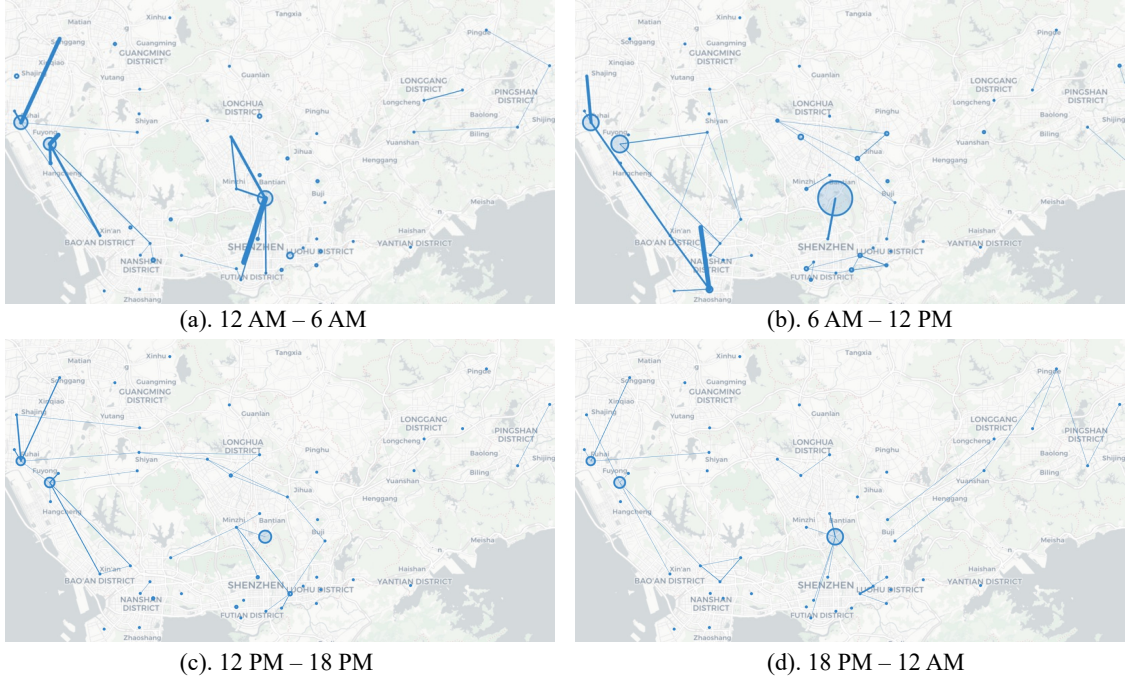

(a). 12 AM – 6 AM  (b). 6 AM – 12 PM

(c). 12 PM – 18 PM  (d). 18 PM – 12 AM

**Figure 4:** Workers in hubs and flow across hubs

## 3.3 Costs

We consider four kinds of worker payments in the experiments, including hiring payment, hour payment, waiting payment at hubs and moving payment across hubs. We also take two kinds of penalties into consideration, including parcels' lateness penalty as well as workforce emergency hiring penalty. The unit costs can be found in the Table 1. To be specific, moving payment across hubs depends on the travelling distance between hubs. Workers are paid 10 or 20 Yuan if the moving distance is respectively less than or more than 3000 meters. Similarly, emergency hiring penalty depends on the time difference between assigning time and shift starting time. If a worker is notified of assignment less than 1 hour before the shift starting time, extra 20 Yuan is paid as the emergency hiring penalty. Also, the emergency hiring penalty is 15, 10 and 5 Yuan if the preparation time difference allowed is less than 2, 4 or 8 hours accordingly.

We contrast the cost performance of three scenarios:
- Scenario 1: workers are allowed to move across hubs; shifts are assigned with a rolling horizon.
- Scenario 2: workers are not allowed to move across hubs; shifts are assigned with a rolling horizon.
- Scenario 3: workers are allowed to move across hubs; shifts are assigned at the beginning of a day.

The total costs to cover the workload of a whole day in the three scenarios are summarized in **Table 1**. Allowing workers to travel across hubs saves us 15,762 Yuan. The reason is that by combing short shifts into a long shift, cost of hiring a new worker can be saved even with the extra transportation cost. In addition, 52,475 Yuan is saved by assigning shifts with a rolling horizon. It is because demand uncertainties are not considered when assigning the shifts at the beginning of a day, in which case many parcels are not processed during their maximum dwell time and thus lateness penalty needs to be paid. For future research, long-term optimized planning with short-term heuristic adjustments could be implemented to make improvements from shifts assigned at the beginning of a day.

Table 1: Costs Comparison among the Three Scenarios.

| Cost Type | Unit Price (Yuan) | Cost of Each Scenario (Yuan) | | |
|---|---|---|---|---|
| | | Scenario 1 | Scenario 2 | Scenario 3 |
| Hiring payment | 50/person/day | 274400 | 302000 | 160450 |
| Hourly payment | 20/person/hour | 343180 | 343180 | 330540 |
| Waiting payment at hubs | 5/person/hour | 13080 | 10815 | 9650 |
| Moving payment across hubs | 10 - 20/person | 10380 | 0 | 5710 |
| Lateness penalty | 5/parcel | 0 | 0 | 199240 |
| Emergency hiring penalty | 5 - 20/person | 19630 | 20440 | 7555 |
| Total cost | - | 660670 | 676435 | 713145 |

## 4. Conclusions and Future Research

In this paper, we propose a novel reactive scheduling and relocation model with rolling horizons for efficient and reliable workforce management in the hyperconnected logistic hubs. First, we show that adding moving workers across multiple nearby hubs is possible and cost-saving, especially for logistic networks in urban areas. Second, the proposed rolling horizon approach with dynamically updated predictions is shown helpful to address the presence of workload uncertainty. Future research involves considerations of different predictive scenarios for dynamic decision-making. Robust scheduling with stochastic optimization models can be included for combining long-term optimized planning with short-term heuristic adjustments, to ensure overall good performance and quick response to dynamic events. Furthermore, studies of sequential decision making can be conducted to make reliable and risk-informed decisions under stochastic environment.